\documentclass[11pt, twoside]{article}
\usepackage{amssymb}
\usepackage{amsmath, amssymb}\usepackage{cite}\usepackage{mathrsfs}\usepackage[amsmath, thmmarks]{ntheorem}
\usepackage{listings}
\usepackage[titletoc]{appendix}\allowdisplaybreaks
\usepackage{multirow}

\makeatletter
\renewcommand\theequation{\thesection.\@arabic\c@equation}
\makeatother
\newtheorem{thm}{Theorem}[section]%
\newtheorem{lem}[thm]{Lemma}%
\newtheorem{cor}[thm]{Corollary}%

\newtheorem{Remark}[thm]{Remark}%
\newtheorem{Con}[thm]{Conjecture}%
\newtheorem{Fac}[thm]{Fact}%
\input cyracc.def

\def\f{\noindent}

\def\demo{\f{\bf Proof}\hskip10pt}

\def\qed{\hfill $\Box$}

\begin{document}

\title{{\bf On two conjectures about pattern avoidance of cyclic permutations}}
\footnotetext{E-mail addresses: Junyao$_{-}$Pan@126.com}

\author{{\bf Junyao Pan}\\
{\footnotesize Wuxi University, Wuxi, Jiangsu, 214105, P. R. China}}

\date{}
\maketitle

%\title{\textbf{On the permutations that strongly avoid the pattern $312$ or $231$}}
%\author{Junyao Pan
% \\\\
%School of Sciences, Wuxi University, Wuxi, Jiangsu, 214105, People's Republic of China\\}
%\date {} \maketitle
%
%
%\baselineskip=16pt
%
%
%\vskip0.5cm

\noindent{\small {\bf Abstract:} Let $\pi$ be a cyclic permutation that can be expressed in its one-line form as $\pi = \pi_1\pi_2 \cdot\cdot\cdot \pi_n$ and in its standard cycle form as $\pi = (c_1,c_2, ..., c_n)$ where $c_1=1$. Archer et al. introduced the notion of pattern avoidance of one-line and the standard cycle form for a cyclic permutation $\pi$, defined as both $\pi_1\pi_2 \cdot\cdot\cdot \pi_n$ and its standard cycle form $c_1c_2\cdot\cdot\cdot c_{n}$ avoiding a given pattern. Let $\mathcal{A}_n(\sigma_1,...,\sigma_k; \tau)$ denote the set of cyclic permutations in the symmetric group $S_n$ that avoid each pattern of $\{\sigma_1,...,\sigma_k\}$ in their one-line forms and avoid $\tau$ in their standard cycle forms. In this paper, we obtain some results about the cyclic permutations avoiding patterns in both one-line and cycle forms. In particular, we resolve two conjectures of Archer et al.

\vskip0.2cm
\noindent{\small {\bf Keywords:} Cyclic permutation; Standard cycle form; Pattern avoidance.


\section {Introduction}

Let $S_n$ denote the symmetric group on $[n]=\{1, 2, \ldots , n\}$. It is well-known that every permutation $\pi$ in $S_n$
can be written either in its cycle form as a product of disjoint cycles or in its one-line notation as $\pi = \pi_1\pi_2\cdots \pi_n$, where $\pi_i = \pi(i)$ for all $i \in [n]$. If $\pi$ is composed of a single $n$-cycle, then $\pi$ is called a \emph{cyclic permutation}. Let $\pi= \pi_1\pi_2\cdots \pi_n\in S_n$
and $\tau= \tau_1\tau_2\cdots \tau_k\in S_k$ with $k\leq n$. If there exists a subset of indices $1\leq i_1<i_2<\cdot\cdot\cdot<i_k\leq n$
such that $\pi_{i_s}>\pi_{i_t}$ if and only if $\tau_s>\tau_t$ for all $1\leq s<t\leq k$, then we say that $\tau$ is \emph{contained} in $\pi$
and the subsequence $\pi_{i_1}\pi_{i_2}\cdot\cdot\cdot \pi_{i_k}$ is called an \emph{occurrence} of $\tau$ in $\pi$ and denoted by $\tau\leq \pi$.
For example, $132\leq24153$, because $2,~5,~3$ appear in the same order of size as the letters in $132$. Indeed, the theory of pattern avoidance in permutations was introduced by Knuth in \cite{K}, which has been widely studied for half a century, refer to \cite{B, V}. A lot of attention has been given to the concept of pattern avoidance over the years. Some interesting and relevant results regarding pattern avoidance can be found in \cite{AF,AE, AG, BC, BS, BD,BDJ,BL, Huang, K, P0, P, P1, SS,ZZY}.

Let $\pi$ be a cyclic permutation in $S_n$. It is well-known that $\pi$ can be written either in its one-line notation or in its cycle form, as $\pi = \pi_1\pi_2 \cdot\cdot\cdot \pi_n$ and $\pi = (c_1,c_2, ..., c_n)$ respectively. Note that $\pi$ has \emph{rotational symmetry}, that is, $\pi = (c_i,c_{i+1}, ..., c_n,c_1,c_2,...,c_{i-1})$ for each $1\leq i\leq n$. If $c_i=1$ then we call $(c_i,c_{i+1}, ..., c_n,c_1,c_2,...,c_{i-1})$ is the \emph{standard cycle form} of $\pi$. Archer et al. \cite{AB} introduced the notion of pattern avoidance of \emph{one-line and cycle form} for a cyclic permutation, that is, $\pi$ avoids $\sigma$ in its one-line form and avoids $\tau$ in its standard cycle form if $\pi_1\pi_2 \cdot\cdot\cdot \pi_n$ avoids $\sigma$ and $c_ic_{i+1}\cdot\cdot\cdot c_nc_1c_2\cdot\cdot\cdot c_{i-1}$ avoids $\tau$; Archer et al. \cite{AB1} defined the notion of pattern avoidance of \emph{one-line and all cycle forms} for a cyclic permutation, namely, $\pi$ avoids $\sigma$ in its one-line form and avoids $\tau$ in its all cycle forms if $\pi_1\pi_2 \cdot\cdot\cdot \pi_n$ avoids $\sigma$ and $c_ic_{i+1}\cdot\cdot\cdot c_nc_1c_2\cdot\cdot\cdot c_{i-1}$ avoids $\tau$ for each $1\leq i\leq n$. By rotational symmetry, we note that all cycle forms avoid a pattern is equivalent to the standard cycle form avoids all cyclic shifts of the pattern. Let $a_n(\sigma_1,...,\sigma_k; \tau)=|\mathcal{A}_n(\sigma_1,...,\sigma_k; \tau)|$, where $\mathcal{A}_n(\sigma_1,...,\sigma_k; \tau)$ denotes the set of cyclic permutations in the symmetric group $S_n$ that avoid each pattern of $\{\sigma_1,...,\sigma_k\}$ in their one-line forms and avoid $\tau$ in their standard cycle forms. Archer et al. \cite{AB} proposed some interesting conjectures about $a_n(\sigma_1,...,\sigma_k; \tau)$, as follows:

\begin{Con}\label{pan1-1}\normalfont([2, Further directions for research])
\begin{itemize}
\item $a_n(3421,4321; 213)=F_{2n-3}$ for any integer $n\geq2$, where $F_n$ is the $n$-th Fibonacci number;
\item $a_n(1324,1423; 213)=\binom{n}{3}+1$ for any positive integer $n\geq3$.
\end{itemize}
\end{Con}

Motivated by \cite{P0}, we resolve Conjecture\ \ref{pan1-1} by constructing bijections. Specifically, we prove that for the pattern $\tau=\tau_1\cdot\cdot\cdot\tau_{k-2}21$ of length is at least $4$, $a_n(\tau,4321; 213)=F_{2n-3}$ for any integer $n\geq2$, and $a_n(1324,1423; 213)=\binom{n}{3}+1$ for any positive integer $n\geq3$. As corollary, we derive that $\mathcal{A}_n(3421,4321; 213)$ and $a_n(4321; 213)=F_{2n-3}$ for any integer $n\geq2$. Meanwhile, we also obtained the corresponding symmetric results by applying symmetry. Actually, we note that the method of constructing bijections is also suitable to enumerate $\mathcal{A}_n(\tau,4321; 213)$ for some other patterns $\tau$, such as $a_n(4312,4321; 213)=2^{n-2}$ for any integer $n\geq2$, and $a_n(3412,4321; 213)=P_{n-1}$ for any integer $n\geq1$, where $P_{n-1}$ is the $(n-1)$-th Pell number.

\section {Preliminaries}

It is well-known that the \emph{Fibonacci numbers} are defined by $F_0=0$ and $F_1=1$, and the recurrence relation $F_n=F_{n-1}+F_{n-2}$ for $n\geq2$; and the \emph{Pell numbers} are defined by $P_0=0$ and $P_1=1$, and the recurrence relation $P_n=P_{n-1}+P_{n-2}$ for $n\geq2$. Next we state a simple fact that will be used to enumerate $\mathcal{A}_n(\sigma_1,\sigma_2; 213)$ by inductive method.

\begin{Fac}\label{pan2-0}\normalfont
Let $\pi=(1,c_2, ...,c_{r-1},2,c_{r+1},..., c_n)$ be a cyclic permutation in $\mathcal{A}_n(\sigma_1,\sigma_2; 213)$ or $\mathcal{A}_n(\sigma; 213)$ with $n\geq3$. Then $\{c_2, ...,c_{r-1}\}=\{n-r+3,...,n\}$ and $\{c_{r+1},..., c_n\}=\{3,...,n-r+2\}$.
\end{Fac}
Based on the Fact\ \ref{pan2-0}, we define $$\mathcal{A}_n(\sigma_1,\sigma_2; 213)\big|_2^j=\Big\{\pi\in\mathcal{A}_n(\sigma_1,\sigma_2; 213)\Big|\pi = (1,c_2, ...,c_{j-1},2,c_{j+1},..., c_n)\Big\}.$$
Thereby, we have
\begin{equation}\label{eq1}
a_n(\sigma_1,\sigma_2; 213)=\sum_{j=2}^{n}\Big|\mathcal{A}_n(\sigma_1,\sigma_2; 213)\big|_2^j\Big|.
\end{equation}
On the other hand, let $\pi^{-1}$ denote the inverses of $\pi$ in $S_n$. We note that the cycle permutation has inverse symmetry. In other words, if $\pi\in\mathcal{A}_n(\sigma_1,\sigma_2; 213)$, then $\pi^{-1}\in\mathcal{A}_n(\sigma^{-1}_1,\sigma^{-1}_2; 312)$. Thus, we deduce following remark.

\begin{Remark}\label{pan2-1}\normalfont
$a_n(\sigma_1,\sigma_2; 213)=a_n(\sigma^{-1}_1,\sigma^{-1}_2; 312)$ for any positive integer $n$.
\end{Remark}

Finally, we provide two lemmas for constructing bijections that will be used in confirming the Conjecture\ \ref{pan1-1}. Actually, the first one has been appeared in \cite[Lemma 2.2]{P0}.

\begin{lem}\label{pan2-2}\normalfont
Let $\pi=(1,c_2, ...,c_{j-1},2,3,...,m,m+1)$ and $\pi'=(1,c_2-m, ...,c_{j-1}-m)$, where $2<j\leq n$ and $m=n-j+1$. If $\pi=\pi_13\cdot\cdot\cdot (m+1)1\pi_{m+2}\cdot\cdot\cdot\pi_n$, then $$\pi'=(\pi_1-m)(\pi_{m+2}-m)\cdot\cdot\cdot(\pi_{c_{j-1}-1}-m)(\pi_{c_{j-1}}-1)(\pi_{c_{j-1}+1}-m)\cdot\cdot\cdot(\pi_n-m).$$
Conversely, if $\pi'=(\pi_1-m)(\pi_{m+2}-m)\cdot\cdot\cdot(\pi_{c_{j-1}-1}-m)(\pi_{c_{j-1}}-1)(\pi_{c_{j-1}+1}-m)\cdot\cdot\cdot(\pi_n-m)$, then $$\pi=\pi_13\cdot\cdot\cdot (m+1)1\pi_{m+2}\cdot\cdot\cdot\pi_n.$$
\end{lem}
\demo Let $\pi=\pi_13\cdot\cdot\cdot (m+1)1\pi_{m+2}\cdot\cdot\cdot\pi_n$. Note that $\pi'(1)=c_2-m=\pi_1-m$ and $\pi'(c_i-m)=c_{i+1}-m=\pi_{c_i}-m$ for $2\leq i<j-1$, and $\pi'(c_{j-1}-m)=1=\pi_{c_{j-1}}-1$. Thus, $$\pi'=(\pi_1-m)(\pi_{m+2}-m)\cdot\cdot\cdot(\pi_{c_{j-1}-1}-m)(\pi_{c_{j-1}}-1)(\pi_{c_{j-1}+1}-m)\cdot\cdot\cdot(\pi_n-m).$$
Similarly, if $\pi'=(\pi_1-m)(\pi_{m+2}-m)\cdot\cdot\cdot(\pi_{c_{j-1}-1}-m)(\pi_{c_{j-1}}-1)(\pi_{c_{j-1}+1}-m)\cdot\cdot\cdot(\pi_n-m)$, then $$\pi=\pi_13\cdot\cdot\cdot (m+1)1\pi_{m+2}\cdot\cdot\cdot\pi_n,$$as desired.  \qed

As an application of Lemma\ \ref{pan2-2}, we immediately derive the following remark.
\begin{Remark}\label{pan2-2-1}\normalfont
Let $\tau=\tau_1\cdot\cdot\cdot\sigma_{k-2}21$ with $k\geq4$. If $\pi=(1,c_2, ...,c_{j-1},2,3,...,m,m+1)$ and $\pi'=(1,c_2-m, ...,c_{j-1}-m)$ where $2<j\leq n$ and $m=n-j+1$, then $\pi$ avoids $\tau$ and $4321$ if and only if $\pi'$ avoids $\tau$ and $4321$.
\end{Remark}

\begin{lem}\label{pan2-3}\normalfont
Let $\pi=(1,m+2, ...,n,2,c_{r+1},...,c_n)$ and $\pi'=(1,c_{r+1}-1, ...,c_{n}-1)$, where $3\leq r<n$ and $m=n-r+1$. If $\pi=\pi_1\pi_2\cdot\cdot\cdot\pi_{m+1}(m+3)\cdot\cdot\cdot n2$, then $$\pi'=(\pi_2-1)\cdot\cdot\cdot(\pi_{c_n-1}-1)\pi_{c_n}(\pi_{c_n+1}-1)\cdot\cdot\cdot(\pi_{m+1}-1).$$
Conversely, if $\pi'=(\pi_2-1)\cdot\cdot\cdot(\pi_{c_n-1}-1)\pi_{c_n}(\pi_{c_n+1}-1)\cdot\cdot\cdot(\pi_{m+1}-1)$, then $$\pi=\pi_1\pi_2\cdot\cdot\cdot\pi_{m+1}(m+3)\cdot\cdot\cdot n2.$$
\end{lem}
\demo Let $\pi=\pi_1\pi_2\cdot\cdot\cdot\pi_{m+1}(m+3)\cdot\cdot\cdot n2$. We see that $\pi'(1)=c_{r+1}-1=\pi_2-1$ and $\pi'(c_i-1)=c_{i+1}-1=\pi_{c_i}-1=$ for $r+1\leq i<n$, and $\pi'(c_{n}-1)=1=\pi_{c_n}$. Thereby, $$\pi'=(\pi_2-1)\cdot\cdot\cdot(\pi_{c_n-1}-1)\pi_{c_n}(\pi_{c_n+1}-1)\cdot\cdot\cdot(\pi_{m+1}-1).$$
Similarly, if $\pi'=(\pi_2-1)\cdot\cdot\cdot(\pi_{c_n-1}-1)\pi_{c_n}(\pi_{c_n+1}-1)\cdot\cdot\cdot(\pi_{m+1}-1)$, then $$\pi=\pi_1\pi_2\cdot\cdot\cdot\pi_{m+1}(m+3)\cdot\cdot\cdot n2,$$
as desired.  \qed

\begin{cor}\label{pan2-4}\normalfont
Let $\pi=(1,2,c_{3},...,c_n)$ and $\pi'=(1,c_{3}-1, ...,c_{n}-1)$. If $\pi=2\pi_2\cdot\cdot\cdot\pi_{n}$, then $$\pi'=(\pi_2-1)\cdot\cdot\cdot(\pi_{c_n-1}-1)\pi_{c_n}(\pi_{c_n+1}-1)\cdot\cdot\cdot(\pi_{n}-1).$$
Conversely, if $\pi'=(\pi_2-1)\cdot\cdot\cdot(\pi_{c_n-1}-1)\pi_{c_n}(\pi_{c_n+1}-1)\cdot\cdot\cdot(\pi_{n}-1)$, then $$\pi=2\pi_2\cdot\cdot\cdot\pi_{n}.$$
\end{cor}
\demo Proceeding as in the proof of Lemma\ \ref{pan2-3}, we deduce this corollary.  \qed

As an application of Corollary\ \ref{pan2-4}, we immediately derive the following remark.
\begin{Remark}\label{pan2-4-1}\normalfont
Let $\tau=\tau_1\cdot\cdot\cdot\sigma_{k-2}21$ with $k\geq4$. If $\pi=(1,2,c_{3},...,c_n)$ and $\pi'=(1,c_{3}-1, ...,c_{n}-1)$, then $\pi$ avoids $\tau$ and $4321$ if and only if $\pi'$ avoids $\tau$ and $4321$.
\end{Remark}

\section {$\mathcal{A}_n(\tau,4321; 213)$}

 Firstly, we provide a key lemma that will be used to enumerate some $\mathcal{A}_n(\tau,4321; 213)$.

\begin{lem}\label{pan3-1}\normalfont
If a cyclic permutation $\pi=(1,c_2, ...,c_{n})$ is in $\mathcal{A}_n(4321; 213)$ with $n\geq5$ and $c_2\neq2$, then the elements after $2$ appear in increasing order.
\end{lem}
\demo Let $\pi=(1,c_2, ...,c_{r-1},2,c_{r+1},..., c_n)$ with $2<r\leq n-2$. Based on the Fact\ \ref{pan2-0}, we see $\{c_2, ...,c_{r-1}\}=\{n-r+3,...,n\}$ and $\{c_{r+1},..., c_n\}=\{3,...,n-r+2\}$. It suffices to prove that $c_{r+1}=3,c_{r+2}=4,..., c_n=n-r+2$. Setting $\pi=\pi_1\pi_2 \cdot\cdot\cdot \pi_n$. Note that $\pi_1=c_2$ and $\pi_{c_{r-1}}=2$. If $c_{r+1}\neq3$ and $c_{r+i}=3$, then $\pi_2=c_{r+1}$ and $\pi_{c_{r+i-1}}=3$. In this case, $\pi_1\pi_2\pi_{c_{r+i-1}}\pi_{c_{r-1}}$ is an occurrence of $4321$ in $\pi$, a contradiction. Therefore, we derive $c_{r+1}=3$. Similarly, we infer that $c_{r+2}=4,..., c_n=n-r+2$. Thereby, the elements after $2$ appear in increasing order, as desired.  \qed

\begin{lem}\label{pan3-2}\normalfont
Let $n\geq5$ and $\tau=\tau_1\cdot\cdot\cdot\sigma_{k-2}21$ with $k\geq4$. Then for each $3\leq j\leq n$, we have $$\big|\mathcal{A}_n(\tau,4321; 213)\big|_2^j\big|=a_{j-1}(\tau,4321; 213).$$
\end{lem}
\demo Consider $\Big|\mathcal{A}_n(\tau,4321; 213)\big|_2^j\Big|$ for $3\leq j\leq n$. According to Lemma\ \ref{pan3-1}, we see that every $\pi\in\mathcal{A}_n(\tau,4321; 213)\big|_2^j$ can be expressed as $(1,c_2, ...,c_{j-1},2,3,...,n-j+2)$. For convenience, we set $m=n-j+1$. Now we define a mapping $\mathit{g}$ by the rule that $$\mathit{g}:(1,c_2, ...,c_{j-1},2,3,...,m,m+1)\mapsto(1,c_2-m, ...,c_{j-1}-m).$$ It follows from Lemma\ \ref{pan2-2} and Remark\ \ref{pan2-2-1} that $\mathit{g}$ is a bijection from $\mathcal{A}_n(\tau,4321; 213)\big|_2^j$ to $\mathcal{A}_{j-1}(\tau,4321; 213)$. Therefore, $$\Big|\mathcal{A}_n(\tau,4321; 213)\big|_2^j\Big|=\Big|\mathcal{A}_{j-1}(\tau,4321; 213)\Big|=a_{j-1}(\tau,4321; 213),$$ as desired.  \qed

\begin{lem}\label{pan3-3}\normalfont
Let $n\geq5$ and $\tau=\tau_1\cdot\cdot\cdot\sigma_{k-2}21$ with $k\geq4$. Then $$\Big|\mathcal{A}_n(\tau,4321; 213)\big|_2^2\Big|=a_{n-1}(\tau,4321; 213).$$
\end{lem}
\demo For every $(1,2,c_3, ...,c_{n})\in\mathcal{A}_n(\tau,4321; 213)\big|_2^2$, we define a mapping $\mathit{h}$ by the rule that $$\mathit{h}:(1,2,c_3, ...,c_{n})\mapsto(1,c_3-1, ...,c_{n}-1).$$ It follows from Lemma\ \ref{pan2-4} and Remark\ \ref{pan2-4-1} that $\mathit{h}$ is a bijection from $\mathcal{A}_n(\tau,4321; 213)\big|_2^2$ to $\mathcal{A}_{n-1}(\tau,4321; 213)$, and therefore, $$\Big|\mathcal{A}_n(\tau,4321; 213)\big|_2^2\Big|=\Big|\mathcal{A}_{n-1}(\tau,4321; 213)\Big|=a_{n-1}(\tau,4321; 213),$$ as desired.  \qed

\begin{thm}\label{pan3-4}\normalfont
Let $\tau=\tau_1\cdot\cdot\cdot\sigma_{k-2}21$ with $k\geq4$. Then for any integer $n\geq2$, we have $$a_n(\tau,4321; 213)=a_n(\tau^{-1},4321; 312)=F_{2n-3}.$$
\end{thm}
\demo By Remark\ \ref{pan2-1}, it suffices to consider $a_n(\tau,4321; 213)$. Note that there is only one cycle permutation in $S_2$, and there are exactly two cycle permutations in $S_3$. Therefore, $a_2(\tau,4321; 213)=F_{2\cdot2-3}=F_1=1$ and $a_3(\tau,4321; 213)=F_{2\cdot3-3}=F_3=2$. Additionally, one easily checks that $a_4(\tau,4321; 213)=F_{2\cdot4-3}=F_5=5$. Thus, $a_n(\tau,4321; 213)=F_{2n-3}$ holds for $n=2,3,4$. It follows from Equation\ \ref{eq1} and Lemma\ \ref{pan3-2} and Lemma\ \ref{pan3-3} that $$a_n(\tau,4321; 213)=a_{n-1}(\tau,4321; 213)+\sum_{j=2}^{n-1}a_j(\tau,4321; 213).$$ Proof by induction on $n$, we deduce that $$a_n(\tau,4321; 213)=F_{2n-5}+\sum_{j=2}^{n-1}F_{2j-3}=F_{2n-5}+(F_1+F_3+\cdot\cdot\cdot+F_{2n-5}).$$ Since $F_0=0$, we have $$a_n(\tau,4321; 213)=F_{2n-5}+(F_0+F_1+F_3+\cdot\cdot\cdot+F_{2n-5}).$$ By the recurrence relation $F_n=F_{n-1}+F_{n-2}$ for $n\geq2$, we deduce that $$a_n(\tau,4321; 213)=F_{2n-5}+F_{2n-4}=F_{2n-3},$$
as desired.   \qed

It is worth pointing out that the Conjecture\ \ref{pan1-1} about $\mathcal{A}_n(3421,4321; 213)$ is exactly the case where $\tau=3421$. In addition, we see that if $\tau=4321$ then $\mathcal{A}_n(4321; 213)=\mathcal{A}_n(\tau, 4321; 213)$ for each $\tau=\tau_1\cdot\cdot\cdot\tau_{k-2}21$ with $k\geq4$. So we have the following two corollaries.

\begin{cor}\label{pan3-5}\normalfont
$a_n(3421,4321; 213)=a_n(4312,4321; 312)=F_{2n-3}$ for any integer $n\geq2$.
\end{cor}

\begin{cor}\label{pan3-6}\normalfont
$a_n(4321; 213)=a_n(4321; 312)=F_{2n-3}$ for any integer $n\geq2$.
\end{cor}

Finally, we consider $\mathcal{A}_n(\tau,4321; 213)$ for some $\tau\in S_4$.

\begin{thm}\label{pan3-7}\normalfont
$a_n(4312,4321; 213)=a_n(3421,4321; 312)=2^{n-2}$ for any integer $n\geq2$.
\end{thm}
\demo By Remark\ \ref{pan2-1}, it suffices to consider $a_n(4312,4321; 213)$. Clearly, $a_2(4312,4321; 213)=2^{2-2}=1$ and $a_3(4312,4321; 213)=2^{3-2}=2$. One easily checks that $a_4(4312,4321; 213)=2^{4-2}=4$. Thus, $a_n(4312,4321; 213)=2^{n-2}$ holds for $n=2,3,4$. By Lemma\ \ref{pan3-1}, we see that $\mathcal{A}_n(4312,4321; 213)\big|_2^j=\emptyset$ for $2<j<n$. Consider $\mathcal{A}_n(4312,4321; 213)\big|_2^2$. In this case, an argument similar to the used in Lemma\ \ref{pan3-3} shows that $\Big|\mathcal{A}_n(4312,4321; 213)\big|_2^2\Big|=\Big|\mathcal{A}_{n-1}(4312,4321; 213)\Big|$. Similarly, we infer that $\Big|\mathcal{A}_n(4312,4321; 213)\big|_2^n\Big|=\Big|\mathcal{A}_{n-1}(4312,4321; 213)\Big|$. By Equation\ \ref{eq1}, it follows that $$a_n(4312,4321; 213)=2a_{n-1}(4312,4321; 213).$$ Proof by induction on $n$, we deduce that $$a_n(4312,4321; 213)=22^{n-3}=2^{n-2},$$ 
as desired.   \qed

\begin{thm}\label{pan3-8}\normalfont
$a_n(3412,4321; 213)=a_n(3412,4321; 312)=P_{n-1}$ for any integer $n\geq1$.
\end{thm}
\demo By Remark\ \ref{pan2-1}, it suffices to consider $a_n(3412,4321; 213)$. Clearly, $a_2(3412,4321; 213)=P_{2-1}=P_1=1$ and $a_3(3412,4321; 213)=P_{3-1}=P_2=2$. Furthermore, one easily checks that $a_4(3412,4321; 213)=P_{4-1}=P_3=5$. Thus, $a_n(3412,4321; 213)=P_{n-1}$ holds for $n=2,3,4$. By Lemma\ \ref{pan3-1}, we see that $\mathcal{A}_n(3412,4321; 213)\big|_2^j=\emptyset$ for $2<j<n-1$. An argument similar to the used in Lemma\ \ref{pan3-3} shows that $$\Big|\mathcal{A}_n(3412,4321; 213)\big|_2^2\Big|=\Big|\mathcal{A}_{n-1}(3412,4321; 213)\Big|,$$ $$\Big|\mathcal{A}_n(3412,4321; 213)\big|_2^n\Big|=\Big|\mathcal{A}_{n-1}(3412,4321; 213)\Big|$$ and $$\Big|\mathcal{A}_n(3412,4321; 213)\big|_2^{n-1}\Big|=\Big|\mathcal{A}_{n-2}(3412,4321; 213)\Big|.$$ By Equation\ \ref{eq1}, it follows that $$a_n(3412,4321; 213)=2a_{n-1}(3412,4321; 213)+a_{n-2}(3412,4321; 213).$$ By the recurrence relation $P_n=2P_{n-1}+P_{n-2}$ for $n\geq2$, we deduce that $$a_n(4312,4321; 213)=P_{n-1},$$
as desired.   \qed

\section {$\mathcal{A}_n(1324,1423; 213)$}
Note that there are exactly two cycle permutations in $S_3$. Thus, $a_3(1324,1423; 213)=\binom{3}{3}+1=2$. Moreover, one easily checks that $a_4(1324,1423; 213)=\binom{4}{3}+1=5$ and $a_5(1324,1423; 213)=\binom{5}{3}+1=11$. Therefore, we see that $a_n(1324,1423; 213)=\binom{n}{3}+1$ holds for $n=3,4,5$. Next we shall prove this conjecture by induction on $n$. Now we start by the following lemma.

\begin{lem}\label{pan4-1}\normalfont
Let $\pi=(1,c_2, ...,c_{r-1},2,c_{r+1},..., c_n)\in\mathcal{A}_n(1324,1423; 213)$ with $n>5$. If $2<r<n$, then either $c_2< \cdot\cdot\cdot<c_{r-1}=n$ or $c_{r+1}<c_{r+2}<\cdot\cdot\cdot<c_n<c_{r-1}=c_n+1$. In addition, $c_{r+1}<c_{r+2}<\cdot\cdot\cdot<c_n<c_{r-1}=c_n+1$ occurs in $3<r<n$. If $r=n$ then $c_{r-1}=n$ or $c_{r-1}=3$; if $r=3$ then $c_{2}=n$.
\end{lem}
\demo Suppose $2<r<n$. Based on the Fact\ \ref{pan2-0}, we see that $\{c_2, ...,c_{r-1}\}=\{n-r+3,...,n\}$ and $\{c_{r+1},..., c_n\}=\{3,...,n-r+2\}$. If $c_{r-1}=n$ then $c_2< \cdot\cdot\cdot<c_{r-1}=n$ because $\pi$ avoids $213$ in its standard cycle form. Suppose $c_{r-1}\neq n$. Let $\pi=\pi_1\pi_2...\pi_n$. Note that $\pi_{c_n}=1$ and $\pi_{c_{r-1}}=2$. If $c_{r-1}>c_n+1$, then $\pi_{c_n}\pi_{c_n+1}\pi_{c_{r-1}}\pi_{c_{r-1}+1}$ is either an occurrence of $1324$ or an occurrence of $1423$ in $\pi$, a contradiction. Thereby, $c_{r-1}=c_n+1$ and so $c_n=n-r+2$. Since $\pi$ avoids $213$ in its standard cycle form, we have $c_{r+1}<c_{r+2}<\cdot\cdot\cdot<c_n$. In addition, it is clear that $c_{r+1}<c_{r+2}<\cdot\cdot\cdot<c_n<c_{r-1}=c_n+1$ occurs in $3<r<n$. Similarly, we infer that if $r=n$ then $c_{r-1}=n$ or $c_{r-1}=3$, and if $r=3$ then $c_{2}=n$ by Fact\ \ref{pan2-0}. The proof is completed.   \qed

Motivated by Lemma\ \ref{pan4-1}, we introduce the following notations. Let $4\leq r\leq n$. We define $$\mathcal{A}^-_n(1324,1423; 213)\big|_2^r=\Big\{\pi\in\mathcal{A}_n(1324,1423; 213)\big|_2^r\Big|\pi=(1,c_2, ...,c_{r-1},2,c_{r+1},..., c_n),c_{r-1}\neq n\Big\}$$ and
$$\mathcal{A}^+_n(1324,1423; 213)\big|_2^r=\Big\{\pi\in\mathcal{A}_n(1324,1423; 213)\big|_2^r\Big|\pi=(1,c_2, ...,c_{r-1},2,c_{r+1},..., c_n),c_{r-1}=n\Big\}.$$
Clearly,
$$\Big|\mathcal{A}_n(1324,1423; 213)\big|_2^r\Big|=\Big|\mathcal{A}^+_n(1324,1423; 213)\big|_2^r\Big|+\Big|\mathcal{A}^-_n(1324,1423; 213)\big|_2^r\Big|.$$
Thereby,
$$\sum_{r=4}^{n}\Big|\mathcal{A}_n(1324,1423; 213)\big|_2^r\Big|=\sum_{r=4}^{n}\Big|\mathcal{A}^+_n(1324,1423; 213)\big|_2^r\Big|+\sum_{r=4}^{n}\Big|\mathcal{A}^-_n(1324,1423; 213)\big|_2^r\Big|.$$

\begin{lem}\label{pan4-2}\normalfont
Let $4\leq r\leq n$ and $n\geq6$. If $r<n$, then $$\Big|\mathcal{A}^+_n(1324,1423; 213)\big|_2^r\Big|=n-r;$$
and if $r=n$ then $$\Big|\mathcal{A}^+_n(1324,1423; 213)\big|_2^n\Big|=1.$$
\end{lem}
\demo Let $\pi=(1,c_2, ...,c_{r-1},2,c_{r+1},..., c_n)\in\mathcal{A}^+_n(1324,1423; 213)\big|_2^r$ and $m=n-r+1$. Since $c_{r-1}=n$ and $\pi$ avoids $213$ in its standard cycle form, we have $c_2<\cdot\cdot\cdot<c_{r-1}=n$. By Fact\ \ref{pan2-0}, we note that $c_2=m+2,c_3=m+3,...,c_{r-1}=n$. Clearly, if $r=n$, then $$\mathcal{A}^+_n(1324,1423; 213)\big|_2^n=\big\{134...n2\big\}~{\rm{ and~ so}}~\Big|\mathcal{A}^+_n(1324,1423; 213)\big|_2^n\Big|=1.$$
Consider $4\leq r<n$. Define a mapping $\mathit{\eta}$ by the rule that $$\mathit{\eta}:(1,m+2, ...,n,2,c_{r+1},...,c_n)\mapsto(1,c_{r+1}-1, ...,c_{n}-1).$$
It follows from Lemma\ \ref{pan2-3} that $\pi$ avoids both $1324$ and $1423$ if and only $(1,c_{r+1}-1, ...,c_{n}-1)$ avoids $132$. Therefore, $\mathit{\eta}$ is a bijection from $\mathcal{A}^+_n(1324,1423; 213)\big|_2^r$ to $\mathcal{A}_{m}(132; 213)$, and thus $$\Big|\mathcal{A}^-_n(1324,1423; 213)\big|_2^r\Big|=a_m(132; 213).$$
By \cite[Theorem 3.8]{AB}, we deduce that $$\Big|\mathcal{A}^+_n(1324,1423; 213)\big|_2^r\Big|=n-r,$$
as desired.   \qed

So far, we have seen that
\begin{equation}\label{eq3}
\sum_{r=4}^n\Big|\mathcal{A}^+_n(1324,1423; 213)\big|_2^r\Big|=1+\sum_{r=4}^{n-1}(n-r)=\binom{n-3}{2}+1.
\end{equation}

\begin{lem}\label{pan4-3}\normalfont
 $\Big|\mathcal{A}_n(1324,1423; 213)\big|_2^3\Big|=\Big|\mathcal{A}_{n-2}(1324,1423; 213)\Big|$ holds for $n\geq6$.
\end{lem}
\demo Define a mapping $\mathit{\rho}$ by the rule that $$\mathit{\rho}:(1,n,2,c_4, ...,c_{n})\mapsto(1,c_4-1, ...,c_{n}-1).$$
Proceeding as in the proof of Lemma\ \ref{pan3-2}, we infer that $\mathit{\rho}$ is a bijection from $\mathcal{A}_n(1324,1423; 213)\big|_2^3$ to $\mathcal{A}_{n-2}(1324,1423; 213)$. Thereby, $$\Big|\mathcal{A}_n(1324,1423; 213)\big|_2^3\Big|=\Big|\mathcal{A}_{n-2}(1324,1423; 213)\Big|,$$
as desired.   \qed

\begin{lem}\label{pan4-4}\normalfont
 $\Big|\mathcal{A}_n(1324,1423; 213)\big|_2^2\Big|=2n-6$ holds for $n\geq6$.
\end{lem}
\demo For convenience, we set $$\mathcal{A}_n(1324,1423; 213)\big|_2^2\big|_n^r=\Big\{(1,2,c_3,c_4, ...,c_{n})\in\mathcal{A}_n(1324,1423; 213)\big|_2^2\Big|c_r=n\Big\},$$ where $3\leq r\leq n$.
Since $\pi$ avoids $213$ in its standard cycle form, it follows that $$\Big|\mathcal{A}_n(1324,1423; 213)\big|_2^2\big|_n^n\Big|=\Big|\big\{(1,2,3,...,n)\big\}\Big|=1.$$
Moreover, one easily checks that $\pi=(1,2,c_3, ...,c_{n-2},n, c_n)\in\mathcal{A}_n(1324,1423; 213)\big|_2^2\big|_n^{n-1}$ if and only if $c_3<\cdot\cdot\cdot<c_{n-2}$. Thereby,
$$\Big|\mathcal{A}_n(1324,1423; 213)\big|_2^2\big|_n^{n-1}\Big|=n-3.$$

Consider $\mathcal{A}_n(1324,1423; 213)\big|_2^2\big|_n^{3}$. Let $\pi=(1,2,n,c_4, ...,c_n)\in\mathcal{A}_n(1324,1423; 213)\big|_2^2\big|_n^{3}$. Since $\pi$ avoids both $1324$ and $1423$ in its one-line, we have $\pi=2n(n-1)(n-2)\cdot\cdot\cdot(n-j)1(n-j-1)\cdot\cdot\cdot3$ for some $j$. Clearly, $c_4=3$ and so $c_5=n-1$ because $\pi$ avoids $213$ in its standard cycle form. By the same token, $c_6=4$ and $c_7=n-2$. As an analogy, we deduce that
\begin{align*}
c_m&=
  \begin{cases}
    n-k & \text{if } m =3+2k, \\
    k+1 & \text{if } m =2k.
  \end{cases}
\end{align*}
Therefore, $$\Big|\mathcal{A}_n(1324,1423; 213)\big|_2^2\big|_n^{3}\Big|=1.$$

Consider $\mathcal{A}_n(1324,1423; 213)\big|_2^2\big|_n^r$ with $3< r< n-1$. Let $$\pi=(1,2,c_3, ...,c_{r-1},n,c_{r+1},..., c_n)\in\mathcal{A}_n(1324,1423; 213)\big|_2^2\big|_n^r.$$ We claim that $\pi_i=i$ for $i=3,4,...,r-1$. Note that if $c_3\neq3$, then $c_{r+1}=3$, otherwise, $\pi$ contains either $1324$ or $1423$. Assume that $c_{r+1}=3$. Since $\pi$ avoids $213$ in its standard cycle form, we have $\{c_{r+2},...,c_n\}=\{4,...,n-r+2\}$. In this case, one easily checks that $\pi$ contains an occurrence of $1324$ in $\pi$. Thus, $c_3=3$. Similarly, we infer that $\pi_i=i$ for $i=4,...,r-1$. An argument similar to the used in $\mathcal{A}_n(1324,1423; 213)\big|_2^2\big|_n^{3}$ shows that for $r+1\leq m\leq n$,
\begin{align*}
c_m&=
  \begin{cases}
    r-1+k & \text{if } m =r-1+2k, \\
    n-k & \text{if } m =r+2k.
  \end{cases}
\end{align*}
Thereby, $\Big|\mathcal{A}_n(1324,1423; 213)\big|_2^2\big|_n^r\Big|=1$ for $3< r< n-1$.
Note that $$\Big|\mathcal{A}_n(1324,1423; 213)\big|_2^2\Big|=\sum_{r=3}^{n}\Big|\mathcal{A}_n(1324,1423; 213)\big|_2^2\big|_n^r\Big|=2n-6,$$
as desired.   \qed

\begin{lem}\label{pan4-5}\normalfont
 $\Big|\mathcal{A}^-_n(1324,1423; 213)\big|_2^n\Big|=\Big|\mathcal{A}_{n-1}(1324,1423; 213)\big|_2^{n-1}\Big|$ holds for $n\geq6$.
\end{lem}
\demo Define a mapping $\mathit{\rho}$ by the rule that $$\mathit{\rho}:(1,c_2, ...,c_{n-2},3,2)\mapsto(1,c_2-1, ...,c_{n-2}-1,2).$$
By the case $m=1$ of Lemma\ \ref{pan2-2}, we deduce that $\mathit{\rho}$ is a bijection from $\mathcal{A}^-_n(1324,1423; 213)\big|_2^n$ to $\mathcal{A}_{n-1}(1324,1423; 213)\big|_2^{n-1}$. Thereby, $$\Big|\mathcal{A}^-_n(1324,1423; 213)\big|_2^n\Big|=\Big|\mathcal{A}_{n-1}(1324,1423; 213)\big|_2^{n-1}\Big|,$$
as desired.   \qed

\begin{lem}\label{pan4-6}\normalfont
Let $3<r<n$. Then $\Big|\mathcal{A}^-_n(1324,1423; 213)\big|_2^r\Big|=\Big|\mathcal{A}_{r-1}(1324,1423; 213)\big|_2^{r-1}\Big|$.
\end{lem}
\demo Set $m=n-r+1$. Define a mapping $\mathit{\eta}$ by the rule that $$\mathit{\eta}:(1,c_2, ...,c_{r-2},m+2,2,3,...,m+1)\mapsto(1,c_2-m, ...,c_{r-2}-m,2).$$
Proceeding as in the proof of Lemma\ \ref{pan3-2}, we deduce that $\mathit{\eta}$ is a bijection from $\mathcal{A}^-_n(1324,1423; 213)\big|_2^n$ to $\mathcal{A}_{r-1}(1324,1423; 213)\big|_2^{r-1}$. Thereby, $$\Big|\mathcal{A}^-_n(1324,1423; 213)\big|_2^r\Big|=\Big|\mathcal{A}_{r-1}(1324,1423; 213)\big|_2^{r-1}\Big|,$$
as desired.   \qed

\begin{lem}\label{pan4-7}\normalfont
Then $\Big|\mathcal{A}_n(1324,1423; 213)\big|_2^n\Big|=n-2$ for $n\geq3$.
\end{lem}
\demo Consider $n\geq6$. Note that $$\Big|\mathcal{A}_n(1324,1423; 213)\big|_2^n\Big|=\Big|\mathcal{A}^-_n(1324,1423; 213)\big|_2^n\Big|+\Big|\mathcal{A}^+_n(1324,1423; 213)\big|_2^n\Big|.$$
It follows from Lemma\ \ref{pan4-2} and Lemma\ \ref{pan4-5} that $$\Big|\mathcal{A}_n(1324,1423; 213)\big|_2^n\Big|=\Big|\mathcal{A}_{n-1}(1324,1423; 213)\big|_2^{n-1}\Big|+1.$$
One easily checks that $$\Big|\mathcal{A}_3(1324,1423; 213)\big|_2^3\Big|=1 ~{\rm{and}}~\Big|\mathcal{A}_4(1324,1423; 213)\big|_2^3\Big|=2~{\rm{and}}~\Big|\mathcal{A}_5(1324,1423; 213)\big|_2^3\Big|=3.$$
Therefore,
$$\Big|\mathcal{A}_n(1324,1423; 213)\big|_2^n\Big|=n-2,$$
as desired.   \qed

\begin{thm}\label{pan4-8}\normalfont
$a_n(1324,1423; 213)=a_n(1324,1342; 312)=\binom{n}{3}+1$ for any integer $n\geq3$.
\end{thm}
\demo By Remark\ \ref{pan2-1}, it suffices to prove $a_n(1324,1423; 213)=\binom{n}{3}+1$ for any integer $n\geq3$. Now we prove it by induction on $n$. For convenience, we set $a^+_n=\sum_{r=4}^n\Big|\mathcal{A}^+_n(1324,1423; 213)\big|_2^r\Big|$ and $a^-_n=\sum_{r=4}^n\Big|\mathcal{A}^-_n(1324,1423; 213)\big|_2^r\Big|$. Note that $$a_n(1324,1423; 213)=\Big|\mathcal{A}_n(1324,1423; 213)\big|_2^2\Big|+\Big|\mathcal{A}_n(1324,1423; 213)\big|_2^3\Big|+a^+_n+a^-_n.$$
It follows from Equation\ \ref{eq3} and Lemma\ \ref{pan4-3} and Lemma\ \ref{pan4-4} that
$$a_n(1324,1423; 213)=a_{n-2}(1324,1423; 213)+(2n-6)+\binom{n-3}{2}+1+a^-_n.$$
By Lemma\ \ref{pan4-5} and Lemma\ \ref{pan4-6}, we deduce that $$a^-_n=\sum_{r=3}^{n-1}\Big|\mathcal{A}_r(1324,1423; 213)\big|_2^r\Big|.$$
It follows from Lemma\ \ref{pan4-7} that $$a^-_n=1+2+\cdot\cdot\cdot+(n-3)=\binom{n-2}{2}.$$
According to the inductive assumption, we deduce that $$a_n(1324,1423; 213)=\binom{n-2}{3}+1+(2n-6)+\binom{n-3}{2}+1+\binom{n-2}{2}=\binom{n}{3}+1,$$
as desired.   \qed

\section{Acknowledgement}

We are very grateful to the anonymous referees for their useful suggestions and comments.


\begin{thebibliography}{99}

\bibitem{AF}
N. Alon and E. Friedgut, On the number of permutations avoiding a given pattern, \emph{J. Combin. Theory Ser. A}, 2000,\textbf{89}(1): 133-140.


\bibitem{AB}
K. Archer, E. Borsh, J. Bridges, C. Graves and M. Jeske, Cyclic permutations avoiding patterns in both one-line and cycle forms, \emph{Preprint arXiv:2312.05145}.


\bibitem{AB1}
K. Archer, E. Borsh, J. Bridges, C. Graves and M. Jeske, Pattern-restricted cyclic permutations with a pattern-restricted cycle form, \emph{Enumer. Comb. Appl.}, 2025, \textbf{5}(1), Paper No. S2R3, 14 pp.

\bibitem{AE}
K. Archer and S. Elizalde, Cyclic permutations realized by signed shifts, \emph{J. Comb.}, 2014, \textbf{5}(1): 1-30.

\bibitem{AG}
K. Archer and A. Geary, Powers of permutations that avoid chains of patterns, \emph{Discrete Math.}, 2024, \textbf{347}(9), Paper No. 114040, 14 pp.


\bibitem{BDJ}
J. Baik, P. Deift and K. Johansson, On the distribution of the length of the longest increasing subsequence of random permutations, \emph{J. Amer. Math. Soc.}, 1999, \textbf{12}(4): 1119-1178.

\bibitem{BL}
N. Blitvi\'c and E. Steingr\'imsson, Permutations, Moments, Measures, \emph{Tran. Amer. Math. Soc.}, 2021, \textbf{347}(8): 5473-5508.


\bibitem{B}
M. B\'ona, \emph{Combinatorics of Permutations}, 2nd edition, CRC Press, 2012.

\bibitem{BC}
M. B\'ona and M. Cory, Cyclic permutations avoiding pairs of patterns of length three, \emph{Discrete Math. Theor. Comput. Sci.}, 2019, \textbf{21}(2), Paper No. 8, 15pp.

\bibitem{BS}
M. B\'ona and R. Smith, Pattern avoidance in permutations and their squares, \emph{Discrete Math.}, 2019, \textbf{342}(11): 3194-3200.

\bibitem{BD}
A. Burcroff and C. Defant, Pattern-avoiding permutation powers, \emph{Discrete Math.}, 2020, \textbf{343}(11): 112017.

\bibitem{Huang}
B. Huang, An upper bound on the number of $(132, 213)$-avoiding cyclic permutations, \emph{Discrete Math.}, 2019, \textbf{342}(6): 1762-1771.

\bibitem{K}
D. Knuth, \emph{The Art of Computer Programming}, \textbf{Vol. 1}, Addison-Wesley, 1968.

\bibitem{P0}
J. Y. Pan, On a conjecture about pattern avoidance of cyclic permutations, \emph{B. Malays. Math Sci. So.}, 2025, \textbf{48} (77).

\bibitem{P}
J. Y. Pan, On a conjecture about strong pattern avoidance, \emph{Graphs Combin.}, 2023, \textbf{39}(1), Paper No. 2, 5 pp.

\bibitem{P1}
J. Y. Pan and P. F. Guo, On the permutations that strongly avoid the pattern 312 or 231, \emph{Graphs Combin.}, 2025, \textbf{41}(1), Paper No. 17, 12 pp.


\bibitem{SS}
R. Simion and F. W. Schmidt, Restricted permutations, \emph{European J. Combin.}, 1985, \textbf{6}(4): 383-406.

\bibitem{V}
V. Vatter, \emph{Permutation classes}, in: Mikl\'os B\'ona (Ed.), Handbook of Enumerative Combinatorics, CRC Press, 2015.


\bibitem{ZZY}
Robin D. P. Zhou, Yongchun Zang, and Sherry H. F. Yan, Further refinements of Wilf-equivalence for patterns of length $4$, \emph{J. Combin. Theory Ser. A}, 2024, \textbf{204}, Paper No. 105863, 13 pp.


\end{thebibliography}
\end{document}